\documentclass[11pt,a4paper]{article}
\usepackage[T1]{fontenc}
\usepackage[utf8]{inputenc}
\usepackage{lmodern,amsmath,amssymb,amsthm,mathtools,booktabs,array}
\usepackage[margin=27mm]{geometry}
\usepackage{microtype}
\usepackage[hidelinks]{hyperref}
\usepackage{enumitem}
\setlist[enumerate]{itemsep=3pt,topsep=5pt}
\newtheorem{theorem}{Theorem}[section]
\newtheorem{proposition}[theorem]{Proposition}
\newtheorem{lemma}[theorem]{Lemma}
\newtheorem{corollary}[theorem]{Corollary}
\theoremstyle{definition}
\newtheorem{definition}[theorem]{Definition}
\newtheorem{example}[theorem]{Example}
\newtheorem{remark}[theorem]{Remark}
\DeclareMathOperator{\Der}{Der}
\DeclareMathOperator{\LocDer}{LocDer}
\DeclareMathOperator{\Ann}{Ann}
\DeclareMathOperator{\End}{End}
\DeclareMathOperator{\Hom}{Hom}
\DeclareMathOperator{\RefHull}{Ref}
\newcommand{\K}{\mathbb K}
\newcommand{\C}{\mathbb C}
\newcommand{\spn}[1]{\langle #1\rangle}
\newcommand{\doi}[1]{\href{https://doi.org/#1}{\nolinkurl{#1}}}
\allowdisplaybreaks
\title{Local and 2-Local Derivations of Finite-Dimensional\\Zinbiel Algebras}
\author{%
Arislanbay Abdireymov\thanks{V.~I.~Romanovskiy Institute of Mathematics, Academy of Sciences of Uzbekistan, Tashkent, Uzbekistan. Email: \texttt{aaabdireymov1001@gmail.com}.},
Aloberdi Sattarov\thanks{Department of Business Administration, Turan International University, Namangan, Uzbekistan. Email: \texttt{saloberdi90@mail.ru}.},
Bakhtiyor Yusupov\thanks{V.~I.~Romanovskiy Institute of Mathematics, Academy of Sciences of Uzbekistan, Tashkent, Uzbekistan; Department of Applied Mathematics and Mathematical Physics, Urgench State University, Urgench, Uzbekistan. Email: \texttt{baxtiyor\_yusupov\_93@mail.ru}.}}

\date{}
\begin{document}
\maketitle
\begin{abstract}
We study local and 2-local derivations of finite-dimensional Zinbiel
algebras. Over any field of characteristic zero, every Zinbiel algebra
$A$ with $A^2\ne0$ and $\dim A\le4$ admits a proper local derivation.
The proof uses a filtration argument. We exhibit two five-dimensional
Zinbiel algebras with nonzero multiplication whose derivations are all
central. Every local derivation on either algebra is therefore a
derivation, so the minimum dimension for local rigidity with nonzero
multiplication is five.
Annihilator-valued interpolation gives a general construction of proper
local derivations. When $A^3=0$, the multiplication tensor describes
all derivations and local derivations. For a finite-dimensional
nilpotent algebra $A$ over any field, we show that every 2-local
derivation is a derivation if and only if $\dim(A/A^2)\le1$.
The proof combines the annihilator interpolation of Sh.~Ayupov,
K.~Kudaybergenov and B.~Omirov with a one-generator argument.
When all derivations are central, every 2-local derivation is obtained
by composing a homogeneous map $A/A^2\to\Ann(A)$ with the quotient map,
and every such composition is a 2-local derivation. The five-dimensional
examples extend to locally rigid Zinbiel algebras with nonlinear
2-local derivations in every dimension at least five.
We also show that derivations given by left multiplication form a Lie
ideal of nilpotency class at most two and give sufficient conditions
for the existence of left outer derivations. This ideal is zero in our
locally rigid examples.
\end{abstract}

\section{Introduction}
A Zinbiel algebra is locally rigid if every local derivation is a
derivation. We determine the smallest dimension in which local
rigidity occurs with nonzero multiplication and give a criterion for
$2$-local rigidity in terms of the quotient by the square. The
two-sided annihilator supplies the interpolation used in both
problems. We also study derivations induced by left multiplication.

The theory began with the work of R.~Kadison \cite{kadison} and D.~Larson and A.~Sourour \cite{larson} in operator algebras. The two-point analogue, introduced by P.~\v{S}emrl \cite{semrl}, asks for a common derivation at each pair of elements and does not require the map to be linear. In finite-dimensional Lie theory, semisimple algebras exhibit both local rigidity \cite{ayupovlie} and 2-local rigidity \cite{ayupov2lie}. These results motivate the corresponding questions for other nonassociative algebras, particularly those whose multiplication is nilpotent.

For Leibniz algebras, the relation between local agreement and the multiplication has been examined in several families. Sh.~Ayupov, A.~Khudoyberdiyev and B.~Yusupov \cite{solvable} studied solvable Leibniz algebras with model or abelian nilradicals, establishing rigidity in cases with a maximal-dimensional complementary space and producing proper local derivations in other cases. Sh.~Ayupov, K.~Kudaybergenov and B.~Yusupov \cite{pfiliform} described local derivations of $p$-filiform non-Lie Leibniz algebras and proved the existence of 2-local derivations which are not derivations. Sh.~Ayupov and B.~Yusupov \cite{quasifiliform} determined local derivations and local automorphisms of naturally graded quasi-filiform Leibniz algebras and investigated the associated Lie algebra and Lie group structures. Together, these studies show why the local and two-point problems should be treated separately.

Further comparisons arise in nonassociative and associative algebras. B.~Yusupov and I.~Yuldashev \cite{solvablelie2025} described local derivations on solvable Lie algebras with a filiform nilradical. A.~K.~Alauadinov and B.~B.~Yusupov \cite{schrodinger2026} proved that every local derivation of the space-time Schr\"odinger Lie algebra is a derivation for spatial dimension at least three. In the associative setting, K.~Abdurasulov, Sh.~Ayupov and B.~Yusupov \cite{associative2025} studied local and 2-local derivations of null-filiform, filiform and naturally graded quasi-filiform algebras. Their results include 2-local rigidity in the null-filiform case and proper 2-local derivations in the filiform and quasi-filiform cases. These examples show that the distinction between local and two-point agreement is not specific to one class of multiplications.

Related local-mapping phenomena for $\frac12$-derivations occur in solvable Leibniz algebras \cite{leibnizhalf2025}, naturally graded non-Lie $p$-filiform algebras \cite{pfiliformhalf2025}, solvable Lie algebras \cite{liehalf2026,modelhalf2026}, and filiform-type associative algebras \cite{associativehalf}. These provide comparisons; the present paper concerns ordinary derivations.

A Zinbiel algebra is a bilinear algebra satisfying
\begin{equation}\label{zinbiel}
(xy)z=x(yz)+x(zy).
\end{equation}
These algebras arise in Leibniz cohomology and operadic duality \cite{loday}. A.~Dzhumadil'daev and K.~Tulenbaev \cite{dzhumadil} proved that finite-dimensional complex Zinbiel algebras are nilpotent. D.~Towers \cite{towers} extended this result to arbitrary fields. Thus nilpotency is automatic for the finite-dimensional Zinbiel algebras considered here and need not be imposed as a separate hypothesis. Central extensions and low-dimensional classifications are developed in \cite{centralextensions,classification,degenerations,corrigendum}. For local mappings, S.~Umrzaqov \cite{umrzaqov} investigates local derivations of specific Zinbiel families. F.~Arzikulov, F.~Nabijonova and F.~Urinboyev \cite{arzikulov} establish 2-local rigidity in the null-filiform case and construct nonlinear 2-local derivations in filiform cases. Our arguments use the quotient by the square and the two-sided annihilator to address structural questions beyond these individual families.

The closest structural comparison is the work of Kh.~Atajonov, A.~Khudoyberdiyev and B.~Yusupov \cite{nilpotent2026}. They prove that every finite-dimensional non-abelian nilpotent complex Lie algebra admits a proper local derivation and a proper local automorphism. For non-Lie right Leibniz algebras, these existence results hold in dimension at most seven or nilpotency class at most four; their eight-dimensional example of class five has both local derivation and local automorphism rigidity. The central-image construction in that paper starts with a derivation satisfying $0\ne D(L^2)\subseteq Z(L)$ and composes it with a projection onto $L^2$. In Section~\ref{annlocal} we formulate the same interpolation mechanism for arbitrary bilinear algebras, replacing the two-sided centre by $\Ann(A)$. Its availability is a sufficient condition, and its validity in the Lie setting does not establish it for every Zinbiel algebra. The main distinction here is the sharp Zinbiel threshold: local derivation rigidity with nonzero multiplication already occurs in dimension five. The present paper also addresses the separate 2-local problem and left-multiplication derivations, which are not consequences of the local existence results in \cite{nilpotent2026}.

Two five-dimensional Zinbiel algebras with nonzero square are
locally rigid. In each example the derivations are exactly the linear
maps from the quotient by the square into the annihilator. Their
dimension is optimal.
\begin{theorem}\label{minimum}
Let $\K$ be a field of characteristic zero.
\begin{enumerate}
\item Every finite-dimensional Zinbiel algebra $A$ with $A^2\ne0$ and $\dim A\le4$ admits a proper local derivation.
\item There exists a five-dimensional Zinbiel algebra with $A^2\ne0$ and $\LocDer(A)=\Der(A)$.
\end{enumerate}
Consequently, over every field of characteristic zero, the minimum dimension for local rigidity with nonzero square is five. Nilpotency is automatic by \cite[Theorem~2.4]{towers}.
\end{theorem}
The lower bound is proved directly from the filtration $R_0=A$, $R_{m+1}=AR_m$. Dimension restrictions leave only two possible configurations when $A^3\ne0$. In each configuration we construct a derivation acting nontrivially on $A^2$ with image there in $\Ann(A)$. Composing it with a projection onto $A^2$ gives a proper local derivation. Thus the lower bound does not depend on the completeness of a classification list. The upper bound follows from explicit five-dimensional multiplication tables and complete derivation computations in Section~\ref{sharp}.

The second question concerns maps agreeing with derivations at two points. Sh.~Ayupov, K.~Kudaybergenov and B.~Omirov \cite{ako} already constructed proper 2-local derivations on arbitrary complex algebras with $\dim A/A^2\ge2$ and $\Ann(A)\ne0$. We use that construction over arbitrary fields and combine it with uniqueness on one-generated algebras. Our principal contributions concern the sharp five-dimensional local-rigidity threshold, its classification-free lower bound, and the full description of 2-local maps when all derivations are central.
\begin{theorem}\label{twomain}
Let $A$ be a finite-dimensional nilpotent bilinear algebra over any field. Every 2-local derivation of $A$ is a derivation if and only if $\dim A/A^2\le1$. If $\dim A/A^2\ge2$, $A$ admits a nonadditive 2-local derivation.
\end{theorem}
The proof in Subsection~\ref{twocriterion} applies to arbitrary finite-dimensional nilpotent bilinear algebras. A single generator determines a derivation uniquely, while two independent quotient coordinates allow nonlinear homogeneous interpolation into the annihilator. This is an existence and rigidity criterion, rather than a classification of all 2-local maps on arbitrary algebras. We obtain a full parametrization under the additional hypothesis that every derivation is central: the 2-local maps are exactly the homogeneous maps $A/A^2\to\Ann(A)$, composed with the quotient map. Applied to the five-dimensional examples and their extensions, this shows that $\LocDer=\Der$ can coexist with nonlinear 2-local derivations.

We also distinguish central derivations from derivations induced by left multiplication. The latter form a Lie ideal of nilpotency class at most two. We obtain tests for left outer derivations and show that the locally rigid examples above have no nonzero left inner derivations. Thus existence of left outer derivations and failure of local rigidity are separate questions.

After the preliminaries in Section~\ref{prelim},
Sections~\ref{annlocal} and~\ref{sharp} develop the local derivation
constructions and establish the sharp dimension bound.
Section~\ref{twolocal} treats $2$-local derivations, and
Section~\ref{leftinner} studies left inner and left outer
derivations. The concluding remarks collect the remaining questions.

\section{Preliminaries}\label{prelim}
The quotient by the square and the two-sided annihilator provide
the spaces used for interpolation. We fix the convention for powers
of a bilinear algebra, define local and $2$-local derivations, and
identify the space of central derivations with
$\Hom(A/A^2,\Ann(A))$.

Unless stated otherwise, the base field $\K$ has characteristic zero and all algebras are finite-dimensional. Subsection~\ref{twocriterion} and Theorem~\ref{centralparam} hold over arbitrary fields. All multiplication tables list only the nonzero products.

For a bilinear algebra $A$, let $A^m$ be the span of all products of at least $m$ elements, with arbitrary bracketings. Then
\[
A^1=A,\qquad A^{m+1}\subseteq A^m,\qquad A^iA^j\subseteq A^{i+j}.
\]
This convention agrees with the usual all-bracketing powers: a product of more than $m$ elements can be regarded as a product of $m$ elements by grouping subexpressions. In particular, $A^2=\operatorname{span}\{xy:x,y\in A\}$ and $A^3=A(A^2)+(A^2)A$. Nilpotency means $A^N=0$ for some $N$. We write
\[
\Ann(A)=\{z\in A:zA=Az=0\}.
\]
\begin{definition}
A linear map $D:A\to A$ is a derivation if $D(xy)=D(x)y+xD(y)$ for all $x,y\in A$. A linear map $\Delta:A\to A$ is a local derivation if for every $x\in A$ there is $D_x\in\Der(A)$ satisfying $\Delta(x)=D_x(x)$. A map $\nabla:A\to A$, not necessarily linear, is a 2-local derivation if for every $x,y\in A$ there is $D_{x,y}\in\Der(A)$ such that
\[
\nabla(x)=D_{x,y}(x),\qquad\nabla(y)=D_{x,y}(y).
\]
A local or 2-local derivation is called proper if it is not a derivation.
\end{definition}
\begin{proposition}
The space $\LocDer(A)$ is a vector subspace of $\End(A)$ containing $\Der(A)$. Every local derivation preserves $A^m$ and $\Ann(A)$. Every 2-local derivation is homogeneous over $\K$.
\end{proposition}
\begin{proof}
Linear combinations of pointwise witnesses witness linear combinations of local derivations. Ordinary derivations preserve powers by the Leibniz rule. If $z\in\Ann(A)$, applying that rule to $za=az=0$ shows $D(z)a=aD(z)=0$. Pointwise witnesses transfer these properties to local derivations. For a 2-local derivation, a witness at $(x,tx)$ gives $\nabla(tx)=t\nabla(x)$ for every $t\in\K$.
\end{proof}
For a vector subspace $S\subseteq\End(V)$, its algebraic reflexive hull is
\[
\RefHull(S)=\{T\in\End(V):T(v)\in Sv\text{ for every }v\in V\}.
\]
Thus $\LocDer(A)=\RefHull(\Der(A))$. Checking locality only on basis vectors is insufficient.
\begin{lemma}\label{central}
The space
\[
\mathcal Z(A)=\{T\in\End(A):T(A)\subseteq\Ann(A),\ T(A^2)=0\}
\cong\Hom(A/A^2,\Ann(A))
\]
is contained in $\Der(A)$. If $x\notin A^2$ and $z\in\Ann(A)$, there is $T\in\mathcal Z(A)$ with $T(x)=z$.
\end{lemma}
\begin{proof}
Both sides of the derivation identity vanish for such a map. For the last assertion choose $f\in A^*$ with $f(A^2)=0$ and $f(x)=1$, and set $T(y)=f(y)z$.
\end{proof}
We call the elements of $\mathcal Z(A)$ central derivations. Their relation to left inner derivations is investigated separately in Section~\ref{leftinner}.

\section{Local derivations and the annihilator}\label{annlocal}
The two-sided annihilator gives two criteria for local derivations
of bilinear algebras. A derivation with nonzero annihilator-valued
action on $A^2$ yields a proper local derivation, while centrality
of every derivation implies local rigidity. When $A^3=0$, the
multiplication tensor describes both the derivation and local
derivation spaces.

\subsection{Annihilator-valued interpolation}
The first construction extends the central-image projection argument used for nilpotent Lie and right Leibniz algebras in \cite[Sections~3 and~5]{nilpotent2026} to arbitrary bilinear algebras. It separates the behaviour on $A^2$ from the behaviour outside it. On the square we use a fixed derivation; away from the square, central derivations can be chosen to match any prescribed annihilator-valued image.
\begin{lemma}\label{constructionlemma}
Suppose $D_0\in\Der(A)$ and $\Delta$ is linear, $\Delta(A)\subseteq\Ann(A)$, and $\Delta|_{A^2}=D_0|_{A^2}$. Then $\Delta$ is local. If $\Delta(A^2)\ne0$, it is proper.
\end{lemma}
\begin{proof}
For $x\in A^2$, use $D_0$. For $x\notin A^2$, Lemma~\ref{central} provides a central derivation with value $\Delta(x)$ at $x$. If $\Delta$ were a derivation, its annihilator-valued image would imply $\Delta(xy)=0$ for every product, contrary to $\Delta(A^2)\ne0$.
\end{proof}
\begin{theorem}\label{projection}
Let $A=V\oplus A^2$, and let $P:A\to A^2$ be the associated projection. If $D\in\Der(A)$ satisfies
\begin{equation}\label{nonzeroimage}
0\ne D(A^2)\subseteq\Ann(A),
\end{equation}
then $DP$ is a proper local derivation.
\end{theorem}
\begin{proof}
Apply Lemma~\ref{constructionlemma}: $DP$ has annihilator-valued image and agrees with $D$ on $A^2$.
\end{proof}
\begin{proposition}
Set $H(A)=\{D\in\Der(A):D(A^2)\subseteq\Ann(A)\}$ and $\mathcal R(A)=\{D|_{A^2}:D\in H(A)\}$. Then
\[
\dim\LocDer(A)-\dim\Der(A)\ge\dim\mathcal R(A).
\]
\end{proposition}
\begin{proof}
For $R\in\mathcal R(A)$, the operator $RP$ is local. The map $R\mapsto RP+\Der(A)$ is injective: if $RP$ is a derivation, its annihilator-valued image forces $RP(A^2)=0$, hence $R=0$.
\end{proof}

\begin{theorem}\label{rigidity}
If every derivation of $A$ has image in $\Ann(A)$, then
\[
\LocDer(A)=\Der(A)=\mathcal Z(A)\cong\Hom(A/A^2,\Ann(A)).
\]
\end{theorem}
\begin{proof}
Every annihilator-valued derivation kills $A^2$, so $\Der(A)=\mathcal Z(A)$. Witnesses for a local map $\Delta$ give $\Delta(A)\subseteq\Ann(A)$ and $\Delta(A^2)=0$. Hence $\Delta\in\mathcal Z(A)$.
\end{proof}

\subsection{Algebras with vanishing third power}\label{thirdpower}
When $A^3=0$, the square lies in the annihilator and the multiplication is determined by a bilinear map on a complement of $A^2$. We first construct proper local derivations using this decomposition and then obtain a complete description of the ordinary and local derivation spaces from the multiplication tensor.
\begin{theorem}\label{thirdconstruction}
Suppose $A^3=0\ne A^2$, and write $A=V\oplus W$, where $W=A^2$. The projection $P$ onto $W$ is a proper local derivation. For any $a,b\in\K$, the map
\[
\Delta_{a,b}(v+w)=av+bw
\]
is local and is a derivation exactly when $b=2a$.
\end{theorem}
\begin{proof}
We have $W\subseteq\Ann(A)$. The map $E(v+w)=v+2w$ is a derivation, since its Leibniz identity on $vu$ reads $2vu=vu+vu$. Theorem~\ref{projection} shows that $EP=2P$ is local, and hence so is $P$. Now $\Delta_{a,b}=aE+(b-2a)P$. Since some product is nonzero, its derivation identity is equivalent to $b=2a$.
\end{proof}
Every algebra with $A^3=0$ satisfies \eqref{zinbiel}, so the theorem applies in particular to Zinbiel algebras.

To pass from the preceding construction to all local derivations, we express the derivation identity in block form. The stabilizer of the multiplication tensor controls the two diagonal blocks, while the annihilator allows the off-diagonal block to provide the necessary pointwise interpolation.

Let $\beta:V\times V\to W$ denote multiplication; its values span $W$. Define
\[
\mathfrak g_\beta=\{(S,T)\in\End(V)\oplus\End(W):
T\beta(v,u)=\beta(Sv,u)+\beta(v,Su)\ \forall v,u\in V\}.
\]
Write $\mathfrak g_V$ and $\mathfrak g_W$ for the two coordinate projections of $\mathfrak g_\beta$.
\begin{theorem}
Relative to $A=V\oplus W$, the derivations are exactly
\[
\begin{pmatrix}S&0\\ Q&T\end{pmatrix},\qquad
(S,T)\in\mathfrak g_\beta,\quad Q\in\Hom(V,W).
\]
The local derivations are exactly
\[
\begin{pmatrix}S_0&0\\ Q_0&T_0\end{pmatrix},\qquad
S_0\in\RefHull(\mathfrak g_V),\quad T_0\in\RefHull(\mathfrak g_W),\quad Q_0\in\Hom(V,W).
\]
Consequently,
\begin{align*}
\dim\Der(A)&=\dim\mathfrak g_\beta+\dim V\dim W,\\
\dim\LocDer(A)&=\dim\RefHull(\mathfrak g_V)+\dim\RefHull(\mathfrak g_W)+\dim V\dim W.
\end{align*}
\end{theorem}
\begin{proof}
Preservation of $W$ gives the displayed block form. Since $W\subseteq\Ann(A)$, the derivation identity reduces precisely to the condition on $(S,T)$, while $Q$ is arbitrary. Testing a local operator at $v\in V$ and $w\in W$ gives the necessary reflexive-hull conditions.

Conversely, consider $x=v+w$. If $v=0$, choose $(S,T)\in\mathfrak g_\beta$ with $Tw=T_0w$. If $v\ne0$, choose $(S,T)$ with $Sv=S_0v$, then choose $Q$ so that $Qv=Q_0v+T_0w-Tw$. In each case the resulting derivation agrees with the prescribed operator at $x$. The dimension formulas follow from the independent block parameters.
\end{proof}
\begin{example}
Over $\C$, consider $u_1^2=u_2^2=z$, with all other products zero. Then
\[
\Der(A)=\left\{\begin{pmatrix}a&b&0\\-b&a&0\\p&q&2a\end{pmatrix}\right\},\qquad
\LocDer(A)=\left\{\begin{pmatrix}a&b&0\\-b&a&0\\p&q&t\end{pmatrix}\right\}.
\]
Indeed $\beta(v,u)=v^{\mathsf T}u\,z$, so $S+S^{\mathsf T}=tI_2$. The first projected space is $\{aI_2+bJ\}$ with $J=\begin{pmatrix}0&1\\-1&0\end{pmatrix}$. In an eigenbasis of $J$ this is the full diagonal space, which is reflexive. Thus the dimensions are four and five.
\end{example}

\section{The sharp minimum dimension for local rigidity}\label{sharp}
For Zinbiel algebras with nonzero square over a field of
characteristic zero, the minimum dimension for local rigidity is
five. A filtration gives the lower bound, and four-dimensional
examples illustrate its construction. Two five-dimensional algebras
attain the bound, and their extensions give locally rigid examples
in every larger dimension.

\subsection{The lower bound and four-dimensional examples}
We prove the lower-dimensional assertion in Theorem~\ref{minimum}. The argument uses only the Zinbiel identity, nilpotency and dimension restrictions. The following filtration makes those restrictions explicit and reduces the proof to a small number of configurations.

Set $R_0=A$ and $R_{m+1}=AR_m$. The Zinbiel identity implies
\begin{equation}\label{filtration}
R_{m+1}\subseteq R_m,\qquad R_m\subseteq A^{m+1},\qquad R_iR_j\subseteq R_{i+j+1}.
\end{equation}
For the last inclusion use induction on $i+j$. For $i=0$ it is the definition. If $i>0$, write an element of $R_i$ as a sum of $ax$ with $x\in R_{i-1}$, and use $(ax)y=a(xy+yx)$ and induction. Consequently $A^3=R_1R_0+R_2\subseteq R_2$. If $A$ is nilpotent, each nonzero $R_m$ strictly decreases: equality would persist under left multiplication and contradict eventual vanishing.
\begin{lemma}\label{onefiltration}
If $A=\K e+R_1$ and $e^{[0]}=e$, $e^{[m+1]}=ee^{[m]}$, then $R_m=\K e^{[m]}+R_{m+1}$. In particular, $\dim R_m-\dim R_{m+1}\le1$.
\end{lemma}
\begin{proof}
The assertion for $m=0$ is the hypothesis. Assuming it for $m$, products $(ce+r)(de^{[m]}+q)$, with $r\in R_1$ and $q\in R_{m+1}$, are congruent to $cd e^{[m+1]}$ modulo $R_{m+2}$ by \eqref{filtration}. These products span $R_{m+1}$.
\end{proof}
\begin{theorem}\label{lowerbound}
Every Zinbiel algebra over a field of characteristic zero, of dimension at most four and with $A^2\ne0$, has a proper local derivation.
\end{theorem}
\begin{proof}
By \cite[Theorem~2.4]{towers}, $A$ is nilpotent. Put $r_m=\dim R_m$. If $R_2=0$, then $A^3=0$ and Theorem~\ref{thirdconstruction} applies. Otherwise $4\ge r_0>r_1>r_2>0$.

\emph{Case 1: $r_1=3$.} Then $r_0=4$, and $A=\K e+R_1$. Lemma~\ref{onefiltration} gives the dimensions $(4,3,2,1,0)$. Set $v=ee$, $p=ev$, $h=ep$. Lemma~\ref{onefiltration} shows that their respective classes span $R_1/R_2$, $R_2/R_3$ and $R_3/R_4$. Each quotient is one-dimensional, so these classes are nonzero. Consequently $e,v,p,h$ are linearly independent and form an adapted basis. Since $R_4=0$, \eqref{filtration} gives $AR_3=R_3A=0$, hence $0\ne h\in\Ann(A)$. The identity gives $ve=2p$ and $pe=3h$. Choose $\pi(e)=1$, $\pi(R_1)=0$, and $\sigma(v)=1$, $\sigma(e)=\sigma(R_2)=0$. Then
\[
xp=\pi(x)h,\quad py=3\pi(y)h,\quad \sigma(xy)=\pi(x)\pi(y),\quad\pi(xy)=0.
\]
Thus $D_0(x)=\pi(x)p+4\sigma(x)h$ is a derivation: both sides of its Leibniz identity equal $4\pi(x)\pi(y)h$. Moreover $D_0(A^2)\subseteq\K h$ and $D_0(ee)=4h\ne0$. Apply Theorem~\ref{projection}.

\emph{Case 2: $r_1\le2$.} Now $(r_1,r_2,r_3)=(2,1,0)$. Write $R_1=\K u\oplus\K w$ and $R_2=\K w$. By \eqref{filtration}, $w\in\Ann(A)$ and $R_1u=0$. Choose coordinate functionals $s,t$ with $s(u)=t(w)=1$ and $s(w)=t(u)=0$. Write
\[
xy=\alpha(x,y)u+\beta(x,y)w,\qquad xu=\gamma(x)w,\qquad ux=\delta(x)w.
\]
Here $\alpha$ and $\gamma$ are nonzero and $\gamma(R_1)=0$. The Zinbiel identity is equivalent to
\begin{equation}\label{alphadelta}
\alpha(x,y)\delta(z)=\gamma(x)(\alpha(y,z)+\alpha(z,y)).
\end{equation}

\emph{Case 2a: $\delta\ne0$.} Choose $x_0$ with $g=\gamma(x_0)\ne0$. Symmetry in $y,z$ in \eqref{alphadelta} gives $\alpha(x_0,y)\delta(z)=\alpha(x_0,z)\delta(y)$, so $\alpha(x_0,y)=k\delta(y)$. Here $k\ne0$: otherwise \eqref{alphadelta}, first with $x=x_0$, would imply $\alpha(y,z)+\alpha(z,y)=0$, and then would give $\alpha(x,y)\delta(z)=0$ for all $x,y,z$. Since $\delta\ne0$, this would force $\alpha=0$, a contradiction. Substituting $\alpha(x_0,y)=k\delta(y)$ into \eqref{alphadelta} with $x=x_0$ yields
\[
\alpha(y,z)+\alpha(z,y)=(k/g)\delta(y)\delta(z).
\]
Using this equality in \eqref{alphadelta} and choosing $z$ with $\delta(z)\ne0$ gives $\alpha(x,y)=(k/g)\gamma(x)\delta(y)$. Substitution yields
\[
\delta(y)\delta(z)=\gamma(y)\delta(z)+\gamma(z)\delta(y).
\]
Choosing $z_0$ with $\delta(z_0)\ne0$, first set $y=z=z_0$ and then $z=z_0$ to obtain $\delta=2\gamma$. Thus $\alpha=c\gamma\otimes\gamma$ for some $c\ne0$. The map
\[
D_0(x)=\gamma(x)u+3c^{-1}s(x)w
\]
is a derivation. Indeed $\gamma(xy)=0$ and $s(xy)=c\gamma(x)\gamma(y)$, so
\[
D_0(xy)=3c^{-1}s(xy)w=3\gamma(x)\gamma(y)w.
\]
Its Leibniz right side is $\gamma(x)\delta(y)w+\gamma(y)\gamma(x)w=3\gamma(x)\gamma(y)w$. It satisfies $D_0(A^2)\subseteq\K w$ and $D_0(x_0x_0)=3g^2w\ne0$.

\emph{Case 2b: $\delta=0$.} Equation~\eqref{alphadelta} implies that $\alpha$ is alternating; it vanishes if either argument belongs to $R_1$. Choosing $p,q$ with $\alpha(p,q)\ne0$ shows $p,q,u,w$ are independent, hence $\dim A=4$. Since $\gamma$ induces a nonzero functional on the two-dimensional quotient $A/R_1$, choose a basis $x,y,u,w$ with $\gamma(x)=1$, $\gamma(y)=0$ and $\alpha(x,y)=a\ne0$. The products are
\begin{equation}\label{case2b}
\begin{aligned}
x^2&=b_{11}w,&xy&=au+b_{12}w,&yx&=-au+b_{21}w,\\
y^2&=b_{22}w,&xu&=w.
\end{aligned}
\end{equation}
If $b_{22}\ne0$, set $z=-a/b_{22}$ and define
\[
D_0x=zy-z(b_{12}+b_{21})u,\qquad D_0y=2au,\qquad D_0u=w,\qquad D_0w=0.
\]
On the five pairs $(x,x),(x,y),(y,x),(y,y),(x,u)$, both sides of the rule are respectively $0,aw,-aw,0,0$; on every other basis pair both vanish. Thus $D_0$ is a derivation and $D_0(xy)=aw\ne0$.

If $b_{22}=0$, set
\[
\begin{aligned}
D_0x&=x-b_{11}u,&D_0y&=-y+(b_{12}+b_{21})u,\\
D_0u&=(b_{21}/a)w,&D_0w&=w.
\end{aligned}
\]
On the same five pairs the two sides are $b_{11}w,(b_{12}+b_{21})w,0,0,w$, and both vanish on all remaining pairs. Hence $D_0$ is a derivation and $D_0(xu)=w\ne0$. In both subcases $0\ne D_0(A^2)\subseteq\K w$. Theorem~\ref{projection} finishes all cases.
\end{proof}

\subsubsection{Four-dimensional illustrations}\label{fourdim}
The preceding lower bound is independent of classification. The following four multiplication tables provide concrete illustrations of its constructions and make the pointwise witnesses explicit. In the notation of \cite{centralextensions}, $L_1,L_2$ are $[\mathfrak N_1^{\C}]^1_{01},[\mathfrak N_1^{\C}]^1_{02}$, while $L_3,L_4$ are $[\mathfrak N_1]^1_{01},[\mathfrak N_1]^1_{02}$. All four tables agree in the displayed bases; no basis change is required. Each algebra has basis $e_1,e_2,e_3,e_4$; only the listed products are nonzero:
\begin{align*}
L_1:&\quad e_1^2=e_2,\quad e_1e_2=e_4,\quad e_2e_1=2e_4,\quad e_3^2=e_4;\\
L_2:&\quad e_1^2=e_2,\quad e_1e_2=e_4,\quad e_2e_1=2e_4,\quad e_1e_3=e_4;\\
L_3:&\quad e_1e_2=e_3,\quad e_2e_1=-e_3,\quad e_1e_3=e_4,\quad e_2^2=e_4;\\
L_4:&\quad e_1e_2=e_3,\quad e_2e_1=-e_3,\quad e_1e_3=e_4.
\end{align*}
For all four algebras, $\Ann(L_i)=\K e_4$ and $L_i^4=0\ne L_i^3$. Their squares are $\spn{e_2,e_4}$ for $i=1,2$ and $\spn{e_3,e_4}$ for $i=3,4$.
\begin{proposition}
For each $L_i$, the operator $D_0$ in Table~\ref{examples} is a derivation and the displayed $\Delta$ is a proper local derivation. Unlisted basis images are zero.
\end{proposition}
\begin{table}[htbp]
\centering\small
\renewcommand{\arraystretch}{1.3}
\begin{tabular}{@{}p{12mm}p{55mm}p{28mm}p{37mm}@{}}
\toprule
Algebra & Nonzero images under $D_0$ & Nonzero image under $\Delta$ & Product detecting failure\\
\midrule
$L_1,L_2$ & $D_0e_1=e_2$, $D_0e_2=3e_4$ & $\Delta e_2=3e_4$ & $\Delta(e_1^2)=3e_4$\\
$L_3$ & $D_0e_1=-e_2$, $D_0e_2=2e_3$, $D_0e_3=e_4$ & $\Delta e_3=e_4$ & $\Delta(e_1e_2)=e_4$\\
$L_4$ & $D_0e_1=e_1$, $D_0e_2=-e_2$, $D_0e_4=e_4$ & $\Delta e_4=e_4$ & $\Delta(e_1e_3)=e_4$\\
\bottomrule
\end{tabular}
\caption{Derivations and proper local derivations on four-dimensional examples.}\label{examples}
\end{table}
\begin{proof}
The Zinbiel identity follows by direct substitution in the multiplication tables. For $L_1,L_2$, the nontrivial derivation identity for $D_0$ is
\[
D_0(e_1^2)=3e_4=e_2e_1+e_1e_2;
\]
the rule on the remaining listed products and on zero products is immediate. For $L_3$, the rule on $e_1e_2=e_3$ reads
\[
D_0(e_3)=e_4=(-e_2)e_2+e_1(2e_3)=-e_4+2e_4.
\]
On $e_2e_1=-e_3$ it reads $-e_4=(2e_3)e_1+e_2(-e_2)$. The identities on $e_1e_3=e_4$ and $e_2^2=e_4$ have both sides zero. The only zero product for which cancellation is needed is $e_1^2=0$, whose right side is $(-e_2)e_1+e_1(-e_2)=e_3-e_3=0$. All remaining zero products give zero on both sides. For $L_4$, $D_0$ is diagonal with weights $(1,-1,0,1)$, and each nonzero product has output weight equal to the sum of its input weights.

In every case $\Delta(L_i)\subseteq\Ann(L_i)$ and $\Delta|_{L_i^2}=D_0|_{L_i^2}$. Lemma~\ref{constructionlemma} proves locality. For an explicit witness, let $x=\sum\limits_{j=1}^{4}x_j e_j$ and choose $g$ with $x_g\ne0$, where $g\in\{1,3\}$ for $L_1,L_2$ and $g\in\{1,2\}$ for $L_3,L_4$. Then
\[
D_x(y)=\frac{y_g}{x_g}\Delta(x)
\]
is a central derivation and $D_x(x)=\Delta(x)$. If no such $g$ exists, then $x\in L_i^2$ and one uses $D_x=D_0$. The final column of Table~\ref{examples} shows $\Delta(L_i^2)\ne0$, so $\Delta$ is not a derivation.
\end{proof}

\subsection{Five-dimensional examples}
The upper bound is realized by the algebra $[\mathfrak N_1]^2_{08}$ in \cite{centralextensions}. Its displayed basis and multiplication table agree with \eqref{Atable}, so no change of basis is needed. The multiplication itself is known; here we compute its derivations and establish local rigidity. We verify the Zinbiel identity and compute the full derivation space, so that the sufficient condition in Theorem~\ref{rigidity} can be applied directly.

Let $A$ have basis $e_1,\ldots,e_5$ and products
\begin{equation}\label{Atable}
e_1^2=e_4,\quad e_1e_2=e_3,\quad e_2e_1=-e_3,\quad e_2^2=e_4,\quad e_1e_3=e_5,\quad e_2e_3=e_4.
\end{equation}
\begin{proposition}\label{Aprop}
The algebra $A$ is Zinbiel, and
\[
A^2=\spn{e_3,e_4,e_5},\quad A^3=\spn{e_4,e_5},\quad A^4=0,\quad\Ann(A)=\spn{e_4,e_5}.
\]
Its derivations are precisely
\[
D(e_1)=pe_4+re_5,\qquad D(e_2)=qe_4+se_5,\qquad D(e_3)=D(e_4)=D(e_5)=0.
\]
In particular $\LocDer(A)=\Der(A)$ and $\dim\Der(A)=4$.
\end{proposition}
\begin{proof}
Here $A^2A=0$, and every symmetrized product belongs to $\spn{e_4,e_5}\subseteq\Ann(A)$, so both sides of \eqref{zinbiel} vanish. The powers follow from the table. If $x=\sum\limits_{i=1}^{5}x_i e_i\in\Ann(A)$, then $xe_1=x_1e_4-x_2e_3=0$, giving $x_1=x_2=0$. Now $e_1x=x_3e_5=0$, giving $x_3=0$.

Put $a=D(e_1)$, $b=D(e_2)$, $c=D(e_3)$, $d=D(e_4)$, $f=D(e_5)$, with subscripts denoting coordinates. The rule on the first four products gives
\begin{align*}
d&=2a_1e_4+a_3e_5=(2b_2+b_3)e_4,\\
c&=(a_1+b_2)e_3+(a_2+b_1)e_4+b_3e_5,\\
-c&=-(a_1+b_2)e_3+(a_2+a_3+b_1)e_4.
\end{align*}
Thus $a_3=b_3=0$, $a_1=b_2$, and $a_2+b_1=0$. Hence $c=2a_1e_3$ and $d=2a_1e_4$. The identity on $e_2e_3=e_4$ reads $2a_1e_4=b_1e_5+3a_1e_4$, forcing $a_1=b_1=0$. Then $a_2=b_2=0$ and $c=d=f=0$. The remaining parameters are $a_4,a_5,b_4,b_5$. Conversely, these maps are central derivations by Lemma~\ref{central}. Apply Theorem~\ref{rigidity}.
\end{proof}

For the second five-dimensional example, let $B$ have basis
$e_1,\ldots,e_5$ and products
\begin{equation}\label{Btable}
e_1^2=e_5,\quad e_1e_2=e_3,\quad e_2e_1=-e_3,\quad e_2^2=e_5,\quad e_2e_4=e_5,\quad e_4e_3=e_5.
\end{equation}
\begin{proposition}\label{Bprop}
The algebra $B$ is Zinbiel, with $B^2=\spn{e_3,e_5}$, $B^3=\K e_5$, $B^4=0$ and $\Ann(B)=\K e_5$. Its derivations are precisely
\[
D(e_1)=pe_5,\qquad D(e_2)=qe_5,\qquad D(e_4)=re_5,\qquad D(e_3)=D(e_5)=0.
\]
Thus $\LocDer(B)=\Der(B)$, of dimension three.
\end{proposition}
\begin{proof}
Again $B^2B=0$ and symmetrized products are annihilator-valued, proving \eqref{zinbiel}. The annihilator conditions $xe_1=xe_3=e_4x=0$ give $x_1=x_2=x_4=x_3=0$.

Write $a=D(e_1)$, $b=D(e_2)$, $c=D(e_3)$, $d=D(e_4)$ and $f=D(e_5)$. The relations $e_1^2=e_2^2=e_5$ and $e_1e_2=-e_2e_1=e_3$ imply
\begin{gather*}
f=2a_1e_5=(2b_2+b_4)e_5,\qquad
c=(a_1+b_2)e_3+(a_2+b_1)e_5,\\
2(a_2+b_1)+a_4=0.
\end{gather*}
The zero relations $e_1e_3=e_1e_4=0$ give $a_4=d_2=0$ and $a_2+d_1=0$. Using $e_4e_1=e_4e_2=e_4^2=0$ gives $d_1=a_2=a_3=b_3=d_3=0$. Also $e_2e_3=0$ gives $b_4=0$. Thus $b_1=0$, $b_2=a_1$, $c=2a_1e_3$ and $f=2a_1e_5$. Finally $e_2e_4=e_4e_3=e_5$ gives
\[
2a_1=b_2+d_4=2a_1+d_4,
\]
so $d_4=a_1=b_2=0$. The remaining maps are exactly those displayed and are central. Theorem~\ref{rigidity} completes the proof.
\end{proof}

\begin{proof}[Proof of Theorem~\ref{minimum}]
Theorem~\ref{lowerbound} gives the lower bound. Proposition~\ref{Aprop} gives a five-dimensional algebra with nonzero square and $\LocDer=\Der$; Proposition~\ref{Bprop} provides another. Therefore the least dimension is exactly five. The automatic nilpotency used in the lower bound follows from \cite[Theorem~2.4]{towers}.
\end{proof}
\begin{remark}
Nonzero multiplication is essential in the statement. For a zero-product algebra, $\Der(A)=\End(A)$ and hence $\LocDer(A)=\Der(A)$ in every dimension. Theorem~\ref{minimum} concerns the first occurrence of this equality when $A^2\ne0$.
\end{remark}

\subsection{Locally rigid algebras in higher dimensions}\label{largerdim}
Adding zero-product summands preserves centrality of derivations
and gives locally rigid examples in every dimension at least five.
A separate six-dimensional multiplication table gives another
example with only central derivations.
\begin{proposition}\label{directsum}
If every derivation of $A$ is central and $U$ is a zero-product algebra, then every derivation of $A\oplus U$ is central. Consequently $\LocDer(A\oplus U)=\Der(A\oplus U)$.
\end{proposition}
\begin{proof}
Write $D(a,u)=(D_Aa+Fu,Ga+Hu)$. The derivation identity gives $D_A\in\Der(A)$, $G(A^2)=0$ and $F(U)\subseteq\Ann(A)$, while $H$ is arbitrary. Thus the image lies in $\Ann(A)\oplus U$ and the square is killed. Apply Theorem~\ref{rigidity}.
\end{proof}
\begin{theorem}\label{extension}
For each $n\ge5$, let $A_n=A\oplus\K^{n-5}$ and $B_n=B\oplus\K^{n-5}$, using \eqref{Atable} and \eqref{Btable} and zero multiplication on the added summands. These are Zinbiel algebras of nilpotency index four, with nonzero square and $\LocDer=\Der$. Moreover,
\[
\dim\Der(A_n)=(n-3)^2,\qquad \dim\Der(B_n)=(n-2)(n-4).
\]
\end{theorem}
\begin{proof}
Apply Proposition~\ref{directsum}. For $A_n$, both the quotient by the square and the annihilator have dimension $n-3$. For $B_n$ their dimensions are $n-2$ and $n-4$. Theorem~\ref{rigidity} gives the formulas.
\end{proof}

For a separate six-dimensional example, let $C$ have basis
$e_1,\ldots,e_6$ and products
\begin{equation}\label{Ctable}
\begin{aligned}
e_1^2&=e_6,&e_1e_2&=e_4,&e_1e_3&=e_6,&e_2e_1&=e_5,\\
e_2^2&=e_6,&e_3e_4&=e_6,&e_3e_5&=-e_6.
\end{aligned}
\end{equation}
\begin{proposition}
The algebra $C$ is Zinbiel, $C^2=\spn{e_4,e_5,e_6}$, $C^3=\K e_6$, $C^4=0$, and $\Ann(C)=\spn{e_4+e_5,e_6}$. Its derivations are exactly
\[
D(e_i)=p_i(e_4+e_5)+q_ie_6\quad(1\le i\le3),\qquad D(e_4)=D(e_5)=D(e_6)=0.
\]
Thus $\LocDer(C)=\Der(C)$ and their dimension is six.
\end{proposition}
\begin{proof}
Here $C^2C=0$, while symmetrized products belong to $\spn{e_4+e_5,e_6}\subseteq\Ann(C)$, proving \eqref{zinbiel}. The annihilator is obtained from $xe_2=xe_4=e_3x=0$.

Write $D(e_6)=ke_6$ and $a=D(e_1)$, $b=D(e_2)$, $c=D(e_3)$. The zero relations $e_3e_1=e_3e_2=e_3^2=0$ give $c_2=c_1=0$, $a_4=a_5$, $b_4=b_5$, $c_4=c_5$. The relations $e_1^2=e_2^2=e_1e_3=e_6$ give
\[
a_2=b_1=0,\qquad k=2a_1+a_3=2b_2=a_1+c_3.
\]
Next $D(e_4)=(a_1+b_2)e_4+b_3e_6$ and $D(e_5)=(a_1+b_2)e_5$. The zero products $e_1e_4=e_2e_4=0$ force $a_3=b_3=0$. Hence $a_1=b_2=c_3$ and $k=2a_1$. Finally $e_3e_4=e_6$ gives $k=3a_1$, hence $a_1=0$. The stated maps are precisely the remaining ones, and they are central derivations.
\end{proof}

\section{2-Local derivations}\label{twolocal}
For a finite-dimensional nilpotent algebra $A$, the dimension of
$A/A^2$ determines whether every $2$-local derivation is a derivation.
When all derivations are central, homogeneous maps from this
quotient to the annihilator describe all $2$-local derivations. This
description applies in particular to the locally rigid examples in
Section~\ref{sharp}.

\subsection{A criterion for 2-local rigidity}\label{twocriterion}
Throughout this subsection, $\K$ is an arbitrary field. The existence
argument for proper $2$-local derivations uses the construction of
\cite{ako}. A single generator, when it exists, determines a derivation
uniquely from its value at that generator.
\begin{lemma}
If $A\ne0$ is nilpotent, then $\Ann(A)\ne0$.
\end{lemma}
\begin{proof}
Choose a nonzero product of maximal length. Multiplying it on either side gives products of greater length, which vanish. Thus it belongs to $\Ann(A)$.
\end{proof}
\begin{lemma}\label{onegenerator}
Let $A$ be nilpotent and $A=\K e+A^2$. Then $A$ is generated by $e$. In particular, a derivation which vanishes at $e$ is zero.
\end{lemma}
\begin{proof}
Let $G$ be the subalgebra generated by $e$. We show $A=G+A^m$ for every $m\ge2$. This holds for $m=2$. If it holds for $m$, substitute $x=g+r$ with $g\in G$, $r\in A^m$ in a product of two arbitrary elements. Every term with an $r$ factor belongs to $A^{m+1}$, so $A^2\subseteq G+A^{m+1}$. Since $A=\K e+A^2$, the induction follows. Nilpotency yields $A=G$. A derivation vanishing at $e$ vanishes on every bracketed word in $e$, hence on $A$.
\end{proof}
\begin{theorem}\label{tworigid}
For a finite-dimensional nilpotent bilinear algebra $A$ over any field, the following are equivalent:
\begin{enumerate}
\item every 2-local derivation of $A$ is a derivation;
\item $\dim_{\K}A/A^2\le1$.
\end{enumerate}
If $\dim A/A^2\ge2$, there exists a 2-local derivation which is not additive.
\end{theorem}
\begin{proof}
If $A=A^2$, induction gives $A=A^m$ for every $m$, so $A=0$. Suppose $A\ne0$ and $\dim A/A^2=1$. Choose $e$ with $A=\K e+A^2$. Given a 2-local derivation $\nabla$, choose $D_0$ agreeing with it at $e$. For each $x$, choose a witness $D_{e,x}$. The difference $D_{e,x}-D_0$ vanishes at $e$, hence is zero by Lemma~\ref{onegenerator}. Therefore $\nabla(x)=D_0(x)$ for every $x$.

Now assume $\dim A/A^2\ge2$. Following \cite{ako}, choose independent functionals $\varphi_0,\varphi_1$ vanishing on $A^2$ and $0\ne z\in\Ann(A)$. Define
\begin{equation}\label{nonlinear}
\nabla(x)=g(\varphi_0(x),\varphi_1(x))z,\qquad
g(s,t)=\begin{cases}s^2/t,&t\ne0,\\0,&t=0.\end{cases}
\end{equation}
The function $g$ is homogeneous over every field. For a pair $x,y$, put $p=(\varphi_0(x),\varphi_1(x))$ and $q=(\varphi_0(y),\varphi_1(y))$. If $p,q$ are independent, a linear functional $\ell:\K^2\to\K$ can be prescribed to have values $g(p),g(q)$. If they are dependent, homogeneity ensures that these values are consistent with their linear relation, so such an $\ell$ again exists. Then
\[
D_{x,y}(v)=\ell(\varphi_0(v),\varphi_1(v))z
\]
is a central derivation agreeing with $\nabla$ at both points.

Choose $u_0,u_1$ with $\varphi_i(u_j)=\delta_{ij}$. Then $\nabla(u_0)=\nabla(u_1)=0$ but $\nabla(u_0+u_1)=z\ne0$. This proves nonadditivity and the converse implication.
\end{proof}
\begin{corollary}
Theorem~\ref{tworigid} holds in particular for finite-dimensional Zinbiel algebras. If the quotient by the square has dimension at least two, proper 2-local derivations exist even over $\mathbb F_2$.
\end{corollary}
\begin{remark}
The one-generated argument is also used for null-filiform Leibniz algebras in \cite{ako}. Its conclusion includes the previously studied null-filiform Zinbiel case \cite{arzikulov}, and parallels the null-filiform associative result in \cite{associative2025}. The theorem gives an existence criterion for arbitrary nilpotent algebras; it does not assert that \eqref{nonlinear} exhausts their 2-local derivations.
\end{remark}

\subsection{The case of central derivations}\label{alltwo}
For bilinear algebras over any field satisfying
$\Der(A)=\mathcal Z(A)$, the quotient map and the annihilator give a
parametrization by homogeneous maps.
\begin{theorem}\label{centralparam}
Suppose $\Der(A)=\mathcal Z(A)$. Let $\pi:A\to V=A/A^2$ be the quotient map, and put $Z=\Ann(A)$. The 2-local derivations of $A$ are exactly the maps
\begin{equation}\label{homogeneous}
\nabla=h\circ\pi,\qquad h:V\to Z,\qquad h(tv)=th(v)\quad(t\in\K,\ v\in V).
\end{equation}
Such a map is a derivation if and only if $h$ is linear.
\end{theorem}
\begin{proof}
Let $\nabla$ be 2-local. Its pointwise witnesses give $\nabla(A)\subseteq Z$. If $x-y\in A^2$, a witness at $(x,y)$ kills $x-y$, so $\nabla(x)=\nabla(y)$. Thus $\nabla$ factors through $V$. Its homogeneity implies homogeneity of the induced map $h$.

Conversely, let $h$ be homogeneous and fix $x,y\in A$. Put $v=\pi(x)$ and $w=\pi(y)$. If $v,w$ are independent, prescribe a linear map $T:V\to Z$ on their span by $T(v)=h(v)$ and $T(w)=h(w)$, and extend it to $V$. If they are dependent, the same prescription is consistent by homogeneity; it also holds when one or both are zero. Then $T\pi$ is a central derivation witnessing 2-locality. Finally, a linear $h$ gives a central derivation, and if $h\pi$ is linear, surjectivity of $\pi$ implies that $h$ is linear.
\end{proof}
\begin{remark}
This is an explicit parametrization: choose a nonzero representative $v_L$ of every one-dimensional subspace $L\subseteq V$, choose arbitrary values $z_L\in Z$, and set $h(tv_L)=tz_L$, $h(0)=0$. Every homogeneous map arises this way. Compatibility between distinct lines is required only if the map is to be linear.
\end{remark}
\begin{corollary}
For the algebras $A,B,C$ in \eqref{Atable}, \eqref{Btable}, \eqref{Ctable}, all 2-local derivations have the following forms, where the displayed scalar functions are arbitrary homogeneous functions of degree one:
\begin{align*}
\nabla_A\Bigl(\sum\limits_{i=1}^{5}x_i e_i\Bigr)&=h_1(x_1,x_2)e_4+h_2(x_1,x_2)e_5,\\
\nabla_B\Bigl(\sum\limits_{i=1}^{5}x_i e_i\Bigr)&=h(x_1,x_2,x_4)e_5,\\
\nabla_C\Bigl(\sum\limits_{i=1}^{6}x_i e_i\Bigr)&=k_1(x_1,x_2,x_3)(e_4+e_5)+k_2(x_1,x_2,x_3)e_6.
\end{align*}
They are derivations precisely when all the corresponding scalar functions are linear.
\end{corollary}
\begin{proof}
Apply Theorem~\ref{centralparam} and the computed squares, annihilators and derivation spaces.
\end{proof}
\begin{corollary}
For every $n\ge5$, the algebras $A_n$ and $B_n$ of Theorem~\ref{extension} have $\LocDer=\Der$ and admit nonlinear 2-local derivations. All their 2-local derivations are the homogeneous maps
\[
A_n/A_n^2\to\Ann(A_n),\qquad B_n/B_n^2\to\Ann(B_n)
\]
composed with the quotient maps.
\end{corollary}
\begin{proof}
All derivations are central. The quotient dimensions are $n-3\ge2$ and $n-2\ge3$. Apply Theorems~\ref{tworigid} and~\ref{centralparam}.
\end{proof}

\section{Left inner and left outer derivations}\label{leftinner}
The derivations induced by left multiplication form a Lie ideal
of nilpotency class at most two. We give sufficient conditions for
left outer derivations and compute the ideal for the locally rigid
examples. Throughout this section, $A$ is a finite-dimensional
Zinbiel algebra over a field of characteristic zero.

\subsection{Left multiplication and the Lie ideal of inner derivations}
Left multiplication in a Zinbiel algebra need not satisfy the derivation identity. We therefore restrict the term left inner derivation to those left multiplications which actually are derivations.

For $a,b\in A$, put $L_a(x)=ax$ and $R_b(x)=xb$, and define
\[
\operatorname{InnDer}_L(A)=\Der(A)\cap\{L_a:a\in A\},\qquad
S(A)=\{a\in A:L_a\in\Der(A)\}.
\]
Both are vector spaces, and $\operatorname{InnDer}_L(A)=\{L_a:a\in S(A)\}$. In this section a left outer derivation means a member of $\Der(A)\setminus\operatorname{InnDer}_L(A)$; this terminology refers specifically to the left-multiplication definition. Write
\[
\operatorname{RAnn}(A)=\{z\in A:Az=0\},\qquad
\operatorname{LAnn}(A)=\{z\in A:zA=0\}.
\]

\begin{lemma}\label{leftcriterion}
The operator $L_a$ is a derivation if and only if
\begin{equation}\label{leftidentity}
a(yx)+x(ay)=0\qquad(x,y\in A).
\end{equation}
If $a\in S(A)$, then $L_a(A)\subseteq A^2$, $L_a(\operatorname{RAnn}(A))=0$, and
\begin{equation}\label{leftcomposition}
L_aL_b=-R_{ab}\qquad(b\in A).
\end{equation}
Moreover, $L_a(A^m)\subseteq A^{m+1}$, so $L_a$ is nilpotent when $A$ is nilpotent.
\end{lemma}
\begin{proof}
Using the Zinbiel identity, the defect of the derivation rule is
\[
a(xy)-(ax)y-x(ay)=-a(yx)-x(ay).
\]
This proves the criterion. The image and annihilator assertions follow directly from left multiplication. Taking $y=b$ in \eqref{leftidentity} gives $a(bx)=-x(ab)$, which is \eqref{leftcomposition}. The filtration assertion follows because multiplication by one more factor increases product length.
\end{proof}

\begin{theorem}\label{leftlieideal}
The space $\operatorname{InnDer}_L(A)$ is a Lie ideal of $\Der(A)$, and
\begin{equation}\label{lefttwostep}
[[E_1,E_2],E_3]=0\qquad(E_i\in\operatorname{InnDer}_L(A)).
\end{equation}
Thus its nilpotency class as a Lie algebra is at most two.
\end{theorem}
\begin{proof}
For every derivation $D$ and every $b\in A$, the derivation rule gives
\begin{equation}\label{derleftcommutator}
[D,L_b](x)=D(bx)-bD(x)=D(b)x,
\qquad [D,L_b]=L_{D(b)}.
\end{equation}
If $L_b$ is a derivation, so is its commutator with $D$, proving the ideal assertion.

Write $E_i=L_{a_i}$ with $a_i\in S(A)$. Formula~\eqref{derleftcommutator} yields
\[
[L_{a_i},L_{a_j}]=L_{a_i a_j},\qquad
L_{a_2a_3+a_3a_2}=0.
\]
Since $L_{a_1a_2}$ is a derivation, another application of that formula and the Zinbiel identity gives
\[
[[L_{a_1},L_{a_2}],L_{a_3}]
=L_{(a_1a_2)a_3}
=L_{a_1(a_2a_3+a_3a_2)}
=[L_{a_1},L_{a_2a_3+a_3a_2}]=0.
\]
\end{proof}

\subsection{Central derivations and tests for left outer derivations}
The central maps used to interpolate local and 2-local derivations need not be left outer. The next result identifies exactly when such maps are left inner, and then gives several sufficient tests for left outer derivations.

\begin{proposition}\label{leftcentraltest}
For $a\in A$, $L_a$ is a central derivation if and only if
\[
aA^2=0,\qquad aA\subseteq\Ann(A).
\]
For $J(A)=\{a\in A:aA^2=0,\ aA\subseteq\Ann(A)\}$, the map
\[
J(A)\longrightarrow\mathcal Z(A)\cong\Hom(A/A^2,\Ann(A)),\qquad a\longmapsto L_a
\]
has kernel $\operatorname{LAnn}(A)$. In particular, every central derivation is left inner if and only if this map is surjective. A rank-one central derivation $D(x)=f(x)z$, where $f(A^2)=0$ and $z\in\Ann(A)$, is left inner exactly when the equations $ax=f(x)z$ have a solution $a\in A$ for all $x\in A$.
\end{proposition}
\begin{proof}
Centrality gives the two displayed conditions. Conversely, these conditions make both terms of \eqref{leftidentity} vanish and give the defining properties of $\mathcal Z(A)$. The kernel and the remaining assertions follow from the definition of left multiplication.
\end{proof}

\begin{lemma}\label{leftkillssquare}
If $L_a$ is a derivation and $aA^2=0$, then $aA\subseteq\Ann(A)$.
\end{lemma}
\begin{proof}
For every $x,y$, the Zinbiel identity gives $(ax)y=a(xy+yx)=0$, while \eqref{leftidentity} with $x,y$ interchanged gives $y(ax)=-a(xy)=0$.
\end{proof}

\begin{theorem}\label{outercriteria}
Let $A\ne0$ be nilpotent. Then $\operatorname{InnDer}_L(A)\subsetneq\Der(A)$ whenever at least one of the following conditions holds:
\begin{enumerate}[label=\textup{(\roman*)}]
\item $A^3=0$;
\item $\operatorname{RAnn}(A)\not\subseteq A^2$;
\item some derivation is not a nilpotent operator;
\item some derivation has image not contained in $A^2$;
\item some derivation is nonzero on $\operatorname{RAnn}(A)$;
\item some derivation kills $A^2$ but has image not contained in $\Ann(A)$;
\item there is $b\in A$ with $A(Ab)=0$ but $A(bA)\ne0$;
\item there are $D_1,D_2,D_3\in\Der(A)$ with $[[D_1,D_2],D_3]\ne0$;
\item the map $J(A)\to\mathcal Z(A)$ in Proposition~\ref{leftcentraltest} is not surjective.
\end{enumerate}
\end{theorem}
\begin{proof}
For (i), choose $A=V\oplus A^2$ and set $E(v+w)=v+2w$. The same calculation as in Theorem~\ref{thirdconstruction} shows that $E$ is a derivation, including the zero-product case. Nilpotency and $A\ne0$ imply $A/A^2\ne0$, so $V\ne0$ and $E$ has image outside $A^2$. It is therefore left outer by Lemma~\ref{leftcriterion}.

For (ii), choose $r\in\operatorname{RAnn}(A)\setminus A^2$, a functional $f$ with $f(A^2)=0$, $f(r)=1$, and $0\ne z\in\Ann(A)$. The last choice is possible by nilpotency. The central derivation $D(x)=f(x)z$ satisfies $D(r)=z\ne0$, whereas every left multiplication kills $r$.

Conditions (iii)--(v) follow from Lemma~\ref{leftcriterion}, and (vi) follows from Lemma~\ref{leftkillssquare}. For (vii), direct use of the Zinbiel identity gives
\[
R_b(xy)-R_b(x)y-xR_b(y)=-x(yb).
\]
Thus $R_b$ is a derivation. If $R_b=L_s$ were left inner, \eqref{leftidentity} would give $(yx)b+x(yb)=0$. But $A(Ab)=0$ and the Zinbiel identity reduce this expression to $y(bx)$, contradicting $A(bA)\ne0$. Condition (viii) follows from Theorem~\ref{leftlieideal}, and (ix) follows from Proposition~\ref{leftcentraltest}.
\end{proof}

\begin{example}\label{centralinnerexample}
Let $T$ have basis $u,v$ with $u^2=v$ and all other products zero. Then
\[
\Der(T)=\{D:D(u)=\alpha u+\beta v,\ D(v)=2\alpha v\},\qquad
\operatorname{InnDer}_L(T)=\mathcal Z(T)=\K L_u.
\]
Indeed $L_u(u)=v$, $L_u(v)=0$, and $L_v=0$. Every central derivation is left inner, although the derivation $D(u)=u$, $D(v)=2v$ is left outer. Thus central interpolation alone does not establish existence of left outer derivations in general.
\end{example}

\subsection{Low-dimensional consequences and the locally rigid examples}
The low-dimensional argument already used for local derivations also supplies left outer derivations. The locally rigid examples give a complementary phenomenon: their nonzero derivations are all left outer.

\begin{corollary}\label{lowdimouter}
Every nonzero Zinbiel algebra of dimension at most four over a field of characteristic zero has a left outer derivation.
\end{corollary}
\begin{proof}
If $A^3=0$, use Theorem~\ref{outercriteria}(i). Otherwise use the cases in the proof of Theorem~\ref{lowerbound}. In Case~1, the adapted basis $e,v,p,h$ satisfies $v^2=3h$, since $(ee)v=e(ev+ve)=3h$. Products of length at least five vanish. If its displayed derivation $D_0$ were $L_s$, write $s=s_0e+s_1v+s_2p+s_3h$. The identity $D_0(e)=p$ forces $s_0=0$, $s_1=1/2$, $s_2=0$, so $L_s(v)=3h/2$, contradicting $D_0(v)=4h$.

In Case~2a, $\alpha=c\gamma\otimes\gamma$ and $\delta=2\gamma$. For $s\in S(A)$, apply \eqref{leftidentity} at $x=y=x_0$, where $g=\gamma(x_0)\ne0$. Its left side is $2cg^2\gamma(s)w$, so $\gamma(s)=0$. It follows that $L_s(A)\subseteq\K w$. But the derivation displayed there has $D_0(x_0)=gu+\lambda w$ for some $\lambda\in\K$, outside $\K w$. Hence $D_0$ is left outer.

In Case~2b with $b_{22}\ne0$, the displayed $D_0(x)$ has a nonzero $y$-coordinate and therefore lies outside $A^2=\spn{u,w}$. With $b_{22}=0$, $D_0(w)=w$, so $D_0$ is not nilpotent. Theorem~\ref{outercriteria}(iv) and (iii), respectively, finish these cases.
\end{proof}

\begin{proposition}\label{rigidouterexamples}
For the algebras $A,B,C$ in \eqref{Atable}, \eqref{Btable} and \eqref{Ctable},
\[
\operatorname{InnDer}_L(A)=\operatorname{InnDer}_L(B)
=\operatorname{InnDer}_L(C)=0.
\]
These algebras are locally rigid, and every nonzero derivation is left outer. The same assertions hold for $A_n,B_n$ in Theorem~\ref{extension}.
\end{proposition}
\begin{proof}
All derivations of these algebras are central, so any left inner derivation must have annihilator-valued image and kill the square. For $A$ and $B$, write $a=\sum\limits_{i=1}^{5}a_i e_i$. For $A$, $L_a(e_3)=a_1e_5+a_2e_4$, so killing the square forces $a_1=a_2=0$; all remaining left multiplications are zero. For $B$, the $e_3$-coordinates in $L_a(e_1)=a_1e_5-a_2e_3$ and $L_a(e_2)=a_1e_3+a_2e_5$ force $a_1=a_2=0$. The equation $L_a(e_3)=a_4e_5=0$ then gives $a_4=0$, and again $L_a=0$.

For $C$, write $a=\sum\limits_{i=1}^{6}a_i e_i$. Killing $e_4\in C^2$ gives $a_3=0$. Now $L_a(e_1)=a_1e_6+a_2e_5$ and $L_a(e_2)=a_1e_4+a_2e_6$ must belong to $\spn{e_4+e_5,e_6}$, forcing $a_2=a_1=0$. Thus $L_a=0$. In a direct sum with a zero-product algebra, left multiplication restricts to that of the original summand and vanishes on the added summand. Its restriction must be a derivation, so the preceding conclusions also apply to $A_n,B_n$. Their local rigidity was established above.
\end{proof}

\subsection{Necessary conditions for equality and the remaining question}
The sufficient tests constrain any possible algebra whose derivations are all left inner. They do not prove that such an algebra cannot exist in arbitrary dimension.

\begin{corollary}\label{leftnecessary}
If $A\ne0$ is a Zinbiel algebra and $\Der(A)=\operatorname{InnDer}_L(A)$, then $\dim A\ge5$, $A^3\ne0$, every derivation is nilpotent, every derivation maps $A$ into $A^2$ and kills $\operatorname{RAnn}(A)$, and
\[
\Ann(A)\subseteq\operatorname{RAnn}(A)\subseteq A^2.
\]
The Lie algebra $\Der(A)$ has nilpotency class at most two, every derivation killing $A^2$ is central, and $J(A)\to\mathcal Z(A)$ is surjective.
\end{corollary}
\begin{proof}
Apply Lemmas~\ref{leftcriterion} and~\ref{leftkillssquare}, Theorem~\ref{leftlieideal}, Theorem~\ref{outercriteria} and Corollary~\ref{lowdimouter}.
\end{proof}

\begin{remark}
The assertion that every derivation is nilpotent is commonly called characteristic nilpotency. For a positively graded algebra $A=\bigoplus_{j\ge1}A_j$, the grading operator $E|_{A_j}=j\,\mathrm{id}$ is a nonnilpotent derivation in characteristic zero. Such an algebra cannot satisfy the equality in Corollary~\ref{leftnecessary}.
\end{remark}

The remaining question is whether a nonzero finite-dimensional complex Zinbiel algebra can satisfy $\Der(A)=\operatorname{InnDer}_L(A)$. Equivalently, does every such algebra have a left outer derivation? The results in this section answer the question under the stated sufficient conditions and in dimensions at most four, but do not resolve it in general.

\section{Concluding remarks}
The examples $A,B,C$ combine local rigidity with nonlinear
$2$-local derivations, and all their nonzero derivations are left
outer. The remaining questions concern the classification of
locally rigid algebras, the full description of $2$-local
derivations, and the existence of left outer derivations.

Condition~\eqref{nonzeroimage} is sufficient for the existence of
proper local derivations, and centrality of all derivations is
sufficient for local rigidity. A classification of the
five-dimensional Zinbiel algebras satisfying $\LocDer=\Der$ remains
a separate problem. The examples $A_n,B_n$ are decomposable for $n>5$,
which leads to the question of whether indecomposable locally rigid
examples exist in arbitrarily large dimensions.

Theorem~\ref{centralparam} describes all $2$-local derivations when
every derivation is central. The corresponding description for
other derivation spaces and the question of when $\LocDer(A)$ is
closed under commutators remain to be investigated.

Section~\ref{leftinner} gives conditions for the existence of left
outer derivations. The general question is whether every nonzero
finite-dimensional complex Zinbiel algebra has a left outer
derivation.

\end{document}